\documentclass[12pt]{article}
\usepackage{latexsym, amsfonts, amsmath}

\usepackage{graphicx}

\newtheorem{thm}{Theorem}[section]

\newtheorem{lem}[thm]{Lemma}

\newtheorem{rem}[thm]{Remark}

\def\infint{\int_{-\infty}^\infty}

\def\convd{\stackrel{\cal D}{\rightarrow}}
\def\convp{\stackrel{P}{\rightarrow}}

\def\ex{{\rm E\,}}

\def\MISE{{\rm MISE_n}}

\def\var{\mathop{\rm Var}\nolimits}

\begin{document}

\bibliographystyle{plain}

\title {Combining Kernel Estimators in\\
the Uniform Deconvolution Problem}

\author {Bert van Es\\[0.3cm]
{\normalsize Korteweg-de Vries Institute  for Mathematics}\\
{\normalsize University of Amsterdam}\\
{\normalsize Science Park 904,
 1018 TV Amsterdam,}\\
{\normalsize P.O. Box  94248, 1090 GE Amsterdam}\\
{\normalsize The Netherlands}}


\maketitle

\begin{abstract}
We construct a density estimator and an estimator of the distribution function
in the uniform deconvolution model.
The estimators are based on inversion formulas and  kernel estimators
of the density of the observations and its derivative. Initially the  inversions yield two different  
estimators of the density and two estimators of the distribution function. We construct   asymptotically optimal convex combinations of these two estimators.
We also derive pointwise  asymptotic
normality of the resulting estimators, the pointwise asymptotic biases  and an expansion
of the mean integrated squared error of the density estimator.
It turns out that the pointwise limit distribution of the
density estimator is the same as the pointwise limit distribution  of the density estimator 
introduced by Groeneboom and Jongbloed (2003), a kernel smoothed nonparametric maximum likelihood estimator of the distribution function.
\\[.5cm]

\end{abstract}

\section{Introduction}

Consider the general deconvolution model. Let $X_1,\ldots, X_n$ be i.i.d. observations, where $ X_i=Y_i+Z_i $
and $Y_i$ and $Z_i$ are independent.  Assume that the unobservable
$Y_i$ have distribution function $F$ and density $f$. Also assume that
the unobservable random variables $Z_i$ have a known density $k$. Note that the
density $g$ of $X_i$ is equal to the convolution of $f$ and $k$, so $g=k*f$ where
$*$ denotes convolution.  So we have
\begin{equation}\label{convgen}
g(x)=\int_{-\infty}^\infty k(x-u)f(u)du.
\end{equation}
The deconvolution problem is the problem of estimating $f$ or $F$
from the observations $X_i$. Later on we will restrict ourselves to {\em uniform deconvolution} where
we require the distribution of the $Z_i$ to be uniform.

\bigskip

Several generally applicable methods have been proposed for this deconvolution model but let us review
{\em direct kernel density estimation} first.
Consider estimation of the density function $g$ from the
observations $X_1,\cdots,X_n$.
The {\em kernel density estimator} with   {\em kernel function}  $w$
and {\em bandwidth}  $h>0$, is defined by
\begin{equation}\label{directest}
g_{nh}(x)=\frac{1}{n}\,\sum_{j=1}^n \frac{1}{h}\,
w\Big(\frac{x-X_j}{h}\Big).
\end{equation}
For smooth $g$, essentially twice continuously differentiable, and
symmetric $w$ with integral one, we have
\begin{align*}
&\ex g_{nh}(x)=\int_{-\infty}^\infty \frac{1}{h}\, w\Big(
\frac{x-u}{h}\Big)g(u)du\\
&\quad = g(x) + \frac{1}{2}\, h^2 g''(x)\int u^2w(u)du + o(h^2),\\
&\var g_{nh}(x) = \frac{1}{nh}\, g(x) \int w^2(u)du +
o\Big(\frac{1}{nh}\Big),
\end{align*}
as $n\to \infty, h\to 0$ and $nh\to \infty$.  For more  on direct kernel density estimators
see for instance Prakasa Rao (1983), Silverman (1986) and Wand and Jones (1995).

 \bigskip

The   standard {\em Fourier type kernel density estimator}  for deconvolution
problems is based on the Fourier transform. For an introduction see for instance Wand
and Jones (1995).
Let $w$ denote a {\em kernel function} and $h>0$ a {\em
bandwidth}.
The estimator $f_{nh}(x)$ of the density $f$ at the point $x$ is defined as
$$
f_{nh}(x)=\frac{1}{2\pi} \int_{-\infty}^\infty
e^{-itx} \frac{\phi_w(ht)\phi_{emp}(t)}{ \phi_k(t)}\,dt
=
{1\over nh}\sum_{j=1}^n v_h\Big({{x-X_j}\over h}\Big),
$$
with
$$
v_h(u)={1\over 2\pi}\infint {{\phi_w(s)}\over\phi_k(s/h)}\
e^{-isu}ds,
$$
where
$$
\phi_{emp}(t) = {1\over n}\sum_{j=1}^n e^{itX_j},
$$
and $\phi_w$ and $\phi_k$ denote the characteristic functions of
$w$ and $k$ respectively.

An important condition for these estimators to
be properly defined is that the characteristic function $\phi_k$
of the density $k$ has no zeroes, which renders it useless for
instance for uniform deconvolution. At the same time this shows that uniform deconvolution is a non standard deconvolution problem. In fact, Hu and Ridder (2004)
argue that in economic applications the assumption of no zeros for $\phi_k$ is not
reasonable. If the error distribution  has  bounded support and is symmetric then its
characteristic function will have zeros. They propose an
approximation of the Fourier transform estimator in such cases.
For  other modifications of the Fourier inversion method, applicable to uniform deconvolution, see Hall
and Meister (2007) and Feuerverger, Kim  and  Sun (2008).

 \bigskip

A second general approach is {\em nonparametric maximum likelihood}.
The likelihood of $F$ is equal to
$$
\prod_{j=1}^n g(X_j)=\prod_{j=1}^n \int_{-\infty}^\infty
k(X_j-t)dF(t).
$$
One can try to explicitly determine a distribution function $F$ that maximizes this likelihood. This
works for exponential deconvolution where a unique explicit maximizing $F$ can be derived, see Jongbloed (1998). For a special case of uniform
deconvolution (if the distribution induced by $F$ concentrated on [0,1]) an explicit expression 
for a maximizing $F$ can also be derived. However, from formula (\ref{convun}) below it follows that
the likelihood is determined by the values $F(X_j)-F(X_j-1), i=j,\ldots,n$. 
Hence any $F$ assigning the same probability to the intervals $(X_j-1,X_j]$ will 
have the same likelihood. So here the maximizer is not unique.
In other cases a numerical maximization procedure is
required. For recent results see the thesis of S. Donauer.

\bigskip

 A third general approach is provided by {\em inversion}.
A selected group of deconvolution problems allows explicit   inversion formulas  of
(\ref{convgen}), expressing the density of interest $f$ in terms
 of the density $g$ of the data. In these cases we can
estimate $f$ by substituting for instance a direct   density
estimate of $g$, for instance the kernel estimate (\ref{directest}), in the inversion formula. In Van Es and Kok (1998)
this strategy has been pursued for deconvolution problems where
$k$ equals the exponential density, the Laplace density, and their
repeated convolutions.

\bigskip

In the    uniform deconvolution problem   the error $Z$ is
Uniform$[0,1)$ distributed. So in this particular deconvolution
problem we assume to have i.i.d. observations from the density
\begin{equation}\label{convun}
g(x)=\int_{-\infty}^\infty I_{[0,1)}(x-u)f(u)du =\int_{x-1}^x f(u)du = F(x)-F(x-1).
\end{equation}
Groeneboom and Jongbloed (2003) consider density estimation in this
problem. They propose a kernel density estimator based on the
nonparametric maximum likelihood estimator (NPMLE) of the
distribution function $F$ and derive its asymptotic properties.
Under the restriction that $f$ is concentrated on the interval 
[0,1], and that $f$ is bounded away from zero, its better
performance compared to a more standard kernel estimator, discussed
below, is noted. Our aim is to show that a kernel type estimator
of $f$ can be constructed which, for all $f$, not necessarily  concentrated 
on [0,1], under some smoothness assumptions  has
the same asymptotic bias and variance as the density estimator of Groeneboom and Jongbloed (2003), cf. Theorem \ref{thm:1} and Remark \ref{equal} below.

In this construction an inversion approach is employed. 
The inversion is based on
 (\ref{convun}). In fact this will lead to two inversion formulas,  yielding two possible estimators.
 We will then combine these estimators into an estimator with asymptotically minimal variance.
We also construct an estimator for the distribution function $F$.
The construction is very similar to that of the density estimator.
For other estimators of the distribution function in uniform
deconvolution we refer to Van Es (1991), Groeneboom and Wellner
(1992) and Van Es and Van Zuijlen (1996).

\section{Uniform deconvolution}\label{univariate}

\subsection{Inversion formulas}
Inversion of the relation (\ref{convun})  is relatively simple. Surprisingly we get two different expressions which of course coincide for density functions $g$ of the form (\ref{convun}). 
The   formulas (\ref{Left1}) and (\ref{inversion1}) have already previously been used in Van Es
(1991) and  Groeneboom and Jongbloed (2003).
\begin{lem}\label{Inversion1}
If $g$ is of the form (\ref{convun}) then we have

\begin{align}
& F(x) = \sum_{j=0}^\infty g(x-j) ,\label{Left1}\\
& F(x) = 1 - \sum_{j=1}^\infty g(x+j) .
\end{align}
Furthermore, assuming that $f$ vanishes at plus and minus infinity, and that $g$ is  continuously differentiable, we have
\begin{align}
& f(x) = \sum_{j=0}^\infty g'(x-j),\label{inversion1}\\
& f(x) = - \sum_{j=1}^\infty g'(x+j).\label{inversion2}
\end{align} 
\end{lem}

\noindent{\bf Proof}

Note that formula (\ref{convun}) can be rewritten as $F(x)=g(x) + F(x-1)$ and, replacing $x$ by $x+1$, as $F(x)=-g(x+1)+F(x+1)$.
Iterating these formulas gives the first two inversion formulas for the distribution function $F$. Differentiating these formulas yields the two formulas for the density $f$.
\hfill$\Box$

\subsection{Estimation of the density function}

We construct our estimator using the two inversion formulas of Lemma \ref{Inversion1}.
The fact that the two expressions for $f$ and $F$ in  (\ref{Inversion1})  are equal, if $g$ is of the form (\ref{convun}),
also follows from the fact that
$$
\sum_{j=-\infty}^\infty g(x+j)=\sum_{j=-\infty}^\infty \{F(x+j)-F(x+j-1)\}=1.
$$
 For an arbitrary density $g$, which is not of the form (\ref{convun}),
 the inversions will in general  not yield distribution functions or
densities, nor
will they coincide. In particular, if we substitute a kernel density
estimator like (\ref{directest})
for $g$ then we get different estimators of $f$ from (\ref{inversion1})  and
(\ref{inversion2}). We get
\begin{equation}\label{estimators}
f^{-}_{nh}(x)=\sum_{j=0}^\infty g'_{nh}(x-j)\quad \mbox{and}\quad
f^{+}_{nh}(x)=-\sum_{j=1}^\infty g'_{nh}(x+j).
\end{equation}
The first of these estimators has also been discussed by Groeneboom
and Jongbloed (2003).

\bigskip

We impose the following condition on the kernel function.

\medskip

\noindent{\bf Condition $W_1$}

{\em The function $w$ is a continuously differentiable symmetric
probability density function with support [-1,1].}

\bigskip

\noindent Because of the bounded support of the kernel estimator
$g_{nh}$ the  sums in (\ref{estimators}) are in fact finite sums.
Moreover, $f^{-}_{nh}$ will be periodic with period one for $x$
large enough and $f^{+}_{nh}$ for $x$ small enough. For
instance,   $f^{-}_{nh}$ is equal to the sum of $g_{nh}'(y)$ over
the values $y=x, x-1, x-2, \ldots$. Once $y$ is on the right hand
side of the support of $g_{nh}'$ this sum does not change anymore if
we replace $x$ by $x+1$. Also note that $f^{-}_{nh}$ vanishes for $x$
smaller than the left endpoint of the support of $g_{nh}'$.  Similarly
$f^{+}_{nh}$ vanishes for $x$
larger than right endpoint of the support of $g_{nh}'$.

\bigskip

Let us derive the kernel estimator. Groeneboom and Jongbloed (2003)
show that $f^{-}_{nh}(x)$ is asymptotically normally distributed.
More precisely, as $n\to\infty, h\to 0$ and $nh\to\infty$, they show
$$
\sqrt{nh^3}(f^{-}_{nh}(x)-\ex f^{-}_{nh}(x))\convd
N(0,\sigma_1^2),
$$
with
$$
\sigma_1^2=F(x)\int w'(u)^2du.
$$
However, by a similar proof it follows that
$$
\sqrt{nh^3}(f^{+}_{nh}(x)-\ex f^{+}_{nh}(x))\convd
N(0,\sigma_2^2),
$$
with
$$
\sigma_2^2=(1-F(x))\int w'(u)^2du.
$$
Apparently  the first estimator has a small variance for
small values of $x$ and the second estimator for large values of
$x$. Hence it makes sense to combine the two. Consider
$$
f^{(t)}_{nh}(x) = t f^{-}_{nh}(x) + (1-t)f^{+}_{nh}(x),
$$
for some fixed $0\leq t\leq 1$. The following theorem establishes
asymptotic normality and the asymptotic bias of this estimator. It also contains the results for
the two estimators (\ref{estimators}) above as special cases, taking $t$ equal to zero and one.

\begin{thm}\label{thm:1.1a}
Assume that Condition $W_1$ is satisfied and that $f$ is bounded on a
neighborhood of $x$. Then, as $n\to\infty, h\to 0, nh \to\infty$,
$$
\sqrt{nh^3}(f^{(t)}_{nh}(x)-\ex f^{(t)}_{nh}(x))\convd
N(0,\sigma_t^2)
$$
with
$$
\sigma_t^2= \Big(t^2 F(x) + (1-t)^2(1-F(x))\Big)\int w'(u)^2du.
$$
Furthermore, if $f$ is twice continuously differentiable on a
neighborhood of $x$ then
$$
\ex  f_{nh}(x) = f(x) + \frac{1}{2}h^2f''(x)\int v^2w(v)dv+
o(h^2).
$$\end{thm}

Up to now $t$ has been an arbitrary constant. It turns out that
the expectation of the estimator does not depend on $t$. See (\ref{expt})
in the proof Theorem (\ref{thm:1.1a}).
However, we can minimize the asymptotic variance by choosing a
specific value for $t$. This variance is minimal if $t$ equals
$1-F(x)$. The minimal value is $F(x)(1-F(x))\int w'(u)^2du$.
Furthermore it turns out that
 if we substitute an estimator
$\hat F_n(x)$ of $F$, which we call the {\em initial estimator}, in $1-F(x)$ for $t$, that is consistent in mean squared
error, then we will achieve the minimal variance. So we introduce
\begin{equation}\label{finalest}
f_{nh}(x) =  (1-\hat F_n(x)) f^{-}_{nh}(x) + \hat
F_n(x)f^{+}_{nh}(x).
\end{equation}
The following theorem establishes asymptotic normality and the
asymptotic bias of this estimator. A suitable estimator $\hat
F_n(x)$ will be constructed in the next section.

\begin{thm}\label{thm:1}
Assume that Condition $W_1$ is satisfied, that $f$ is bounded on a
neighborhood of $x$, and that $\hat F_n(x)$ is an estimator of
$F(x)$ with
\begin{equation}\label{fhatcond}
\ex (\hat F_n(x)-F(x))^2\to 0.
\end{equation}
Then, as $n\to\infty, h\to 0, nh \to\infty$, we have
$$
\sqrt{nh^3}(f_{nh}(x)-\ex f_{nh}(x))\convd N(0,\sigma^2),
$$
with
\begin{equation}\label{limitvar}
\sigma^2=F(x)(1-F(x))\int w'(u)^2du.
\end{equation}
Furthermore, if $f$ is twice continuously differentiable on a neighborhood
of $x$ and
\begin{equation}\label{fhatcond2}
\ex (\hat F_n(x)-F(x))^2=o(nh^7)
\end{equation}
then we have
$$
\ex  f_{nh}(x) = f(x) + \frac{1}{2}h^2f''(x)\int v^2w(v)dv+ o(h^2).
$$
\end{thm}

\begin{rem}\label{equal}
The theorem shows that the kernel density estimator $f_{nh}(x)$
has the same asymptotic properties as the density estimator of Groeneboom and Jongbloed (2003) under the restriction that
$f$ is concentrated on the interval [0,1], and that $f$ is bounded
away from zero. However, in Section \ref{GJvar} we show that the limit
variance of the kernel smoothed NPMLE is in fact equal to (\ref{limitvar}),
even if the restriction of the support of $f$ to [0,1] does not hold.
This means that the limit distibutions of the kernel smoothed NPMLE and our estimator coincide.
For the estimators of Hall and Meister (2007) and Feuerverger et al. (2008) the limit
distributions are not known.
\end{rem}

\begin{rem}
Admittedly, the estimator (\ref{finalest}) lacks
the desirable properties that the estimates are nonnegative and
that their integral is equal to one, which are guaranteed for the
kernel smoothed NPMLE.  
\end{rem}

\subsection{Estimation of the distribution function}

To combine the two density estimators in the previous section optimally
we need an estimator  $\hat F_n(x)$ of $F(x)$. The construction of such an estimator
is similar to the construction of the density estimator.

The
inversion formulas (\ref{inversion1}) and (\ref{inversion2}) can
again be used to construct estimators
\begin{equation}\label{Festimators}
F_{nh}^{-}(x) = \sum_{j=0}^\infty g_{nh}(x-j) \quad
\mbox{and}\quad F_{nh}^{+}(x) = 1 - \sum_{j=1}^\infty
g_{nh}(x+j).
\end{equation}
By similar techniques as in the proof of Theorem \ref{thm:1} one can
show
$$
\sqrt{nh}(F^{-}_{nh}(x)-\ex F^{-}_{nh}(x))\convd
N(0,\tau_1^2),
$$
with
$$
\tau_1^2=F(x)\int w(u)^2du.
$$
and
$$
\sqrt{nh}(F^{+}_{nh}(x)-\ex F^{+}_{nh}(x))\convd
N(0,\tau_2^2),
$$
with
$$
\tau_2^2=(1-F(x))\int w(u)^2du.
$$
Now define the estimator $F_{nh}^{(t)}(x)$ by
$$
F_{nh}^{(t)}(x) =  t F^{-}_{nh}(x) + (1-t)F^{+}_{nh}(x).
$$
The following theorem establishes
asymptotic normality and the asymptotic bias of this estimator.

\begin{thm}\label{thm:2.1b} Assume that Condition $W_1$ is satisfied.
Then, as $n\to\infty, h\to 0, nh \to\infty$,
$$
\sqrt{nh}(F^{(t)}_{nh}(x)-\ex F^{(t)}_{nh}(x))\convd N(0,\tau_t^2),
$$
with
$$
\tau_t^2= \Big(t^2 F(x) + (1-t)^2(1-F(x))\Big)\int w(u)^2du.
$$
Furthermore, if $f$ is  continuously differentiable on a
neighborhood of $x$ then
$$
\ex  F_{nh}(x) = F(x) + \frac{1}{2}h^2f'(x)\int v^2w(v)dv+ o(h^2).
$$
\end{thm}

\bigskip

The same steps, i.e. optimizing over $t$, that resulted in the density
estimator (\ref{finalest}) can be repeated to construct an improved
estimator of $F$. Define $F_{nh}(x)$ by
\begin{equation}\label{finalFest}
F_{nh}(x) =  (1-\hat F_n(x)) F^{-}_{nh}(x) + \hat
F_n(x)F^{+}_{nh}(x).
\end{equation}
We get  the following analogue of Theorem \ref{thm:1}. Note that the rate of convergence is faster than in the density estimation case.

\begin{thm}\label{thm:2}
Assume that Condition $W_1$ is satisfied and that $\hat F_n(x)$ is an estimator of $F(x)$ with
$$
\ex (\hat F_n(x)-F(x))^2\to 0.
$$
Then, as $n\to\infty, h\to 0, nh \to\infty$, we have
$$
\sqrt{nh }(F_{nh}(x)-\ex F_{nh}(x))\convd N(0,\tau^2),
$$
with
$$
\tau^2=F(x)(1-F(x))\int w(u)^2du.
$$
Furthermore, if $f$ is  continuously differentiable on a neighborhood of
 $x$ and
\begin{equation}\label{fhatcond4}
\ex (\hat F_n(x)-F(x))^2=o(nh^5)
\end{equation}
then we have
$$
\ex  F_{nh}(x) = F(x) + \frac{1}{2}h^2f'(x)\int v^2w(v)dv+ o(h^2).
$$
\end{thm}

\noindent{\bf Proof}

 The fact that we can replace $t=1-F(x)$ by a
consistent estimator $1-\hat F_n(x)$ and the bias expansion also
follow
 as in the corresponding parts of
the proof of Theorem{ \ref{thm:1}. \hfill$\Box$

\bigskip

As we will see in Section \ref{sim}   the full subtlety of this result is not needed
to combine the two density estimators in the way of the previous section. It turns out that the plain average
$\frac{1}{2}\,(F^{-}_{nh}(x) + F^{+}_{nh}(x))$ suffices for that purpose if we consider pointwise estimation.
For global properties derived  in Section \ref{miseexp} we will see that  $t$ will have to depend on $x$.

\subsection{Mean integrated squared error of the density estimator}\label{miseexp}

Up to now we have considered pointwise, i.e. for a fixed $x$, properties of the estimators.
Assuming that $f$ is square integrable, an important measure of the global performance of a density estimator is the mean integrated squared error, given by
\begin{equation}\label{mise}
\MISE(h) = \ex \int (f_{nh}(x)-f(x))^2dx.
\end{equation}
If we want to consider this global distance then we have to make sure that our estimator $f_{nh}$ is
square integrable. If we use a fixed weight $t$, independent of $x$, for combining $f^{-}_{nh}(x)$ an $f^{+}_{nh}(x)$, then this is certainly not true because of the periodicity
at plus or minus infinity of $f^{-}_{nh}(x)$ and $f^{+}_{nh}(x)$. We can repair this as follows. If we use the optimal
true weight $t=1-F(x)$ the "estimator" is square integrable once $F$ and $1-F$ are square integrable in the left and right tail respectively. This holds because the estimator $f_{nh}^{-}$ has a finite (random) left end point of its support. Similarly
$f^{+}_{nh}(x)$ has a finite right end point of its support. However, we still have to estimate this weight.
As an estimator of the optimal weights we will use an estimator $F_{nh}^{(t)}$ with $t$
dependent on $x$. In particular we will choose $t=1-H(x)$ where $H$ is a distribution function with square integrable tails.
So we will use
\begin{equation}\label{pivot1}
F_{nh}^{(H(x))}(x)=(1-H(x)) F^{-}_{nh}(x) + H(x)F^{+}_{nh}(x).
\end{equation}
By the same reasoning as above $F_{nh}^{(H(x))}$ has square integrable tails, and thus so has the density estimator that uses
this initial estimate of $F$ for the weight. Of course there is no need to use the same bandwidth for the density estimators $f^{-}_{nh}(x)$ and $f^{+}_{nh}(x)$ as for estimating the weights.

The next theorem gives an expansion of (\ref{mise}) which allows us to establish rate optimality of our estimator.

\begin{thm}\label{misethm} Assume that Condition $W_1$ is satisfied,
that $\int F(x)(1-F(x))dx<\infty$, that $f$  is square integrable and twice continuously differentiable
with bounded and square integrable second derivative.
Furthermore, assume that $\hat F_n(x)$ is an estimator of
$F(x)$ with, as $n\to\infty$,
\begin{equation}\label{fhatcond3}
\int(\ex (\hat F_n(x)-F(x))^4)^{1/2}dx
=o(nh^{7}).
\end{equation}
Then, as $n\to\infty, h\to 0$ and $nh\to \infty$ we have
\begin{align}
\MISE(h)&= \frac{1}{4}h^4\int f''(x)^2dx\Big(\int v^2w(v)dv\Big)^2\nonumber +
\frac{1}{nh^3}\int F(x)(1-F(x))dx\int w'(v)^2dv\\
&+o(h^4)+O\Big(\frac{1}{nh^2}\Big).\nonumber
\end{align}
\end{thm}
The next lemma ensures that we can use (\ref{pivot1}) as pivotal estimator.

\begin{lem}\label{pivotok}
Assume that $f$ is differentiable, that $f'$ is bounded and continuous and that $\int \{F(x)(1-F(x))\}^{1/2}dx<\infty$.
Also assume  that $\int f'(x)^2dx<\infty$
and that $\int H(x)^2(1-H(x))^2dx<\infty$, i.e. $H$ and $1-H$ are square integrable in the left and right tail.
Then, if $n\to \infty, h\to 0$ and $nh\to \infty$, we have, with $F_{nh}^{(H(x))}$ equal to the initial estimator (\ref{pivot1}),
\begin{equation}\label{pivotbound}  
\int(\ex (F_{nh}^{(H(x))}(x)-F(x))^4)^{1/2}dx
=O(h^4) + O\Big(\frac{1}{nh}\Big).
\end{equation}
\end{lem}

The lemma shows that if we use a bandwidth of the form $h=c_1n^{-1/5}$ for the initial estimator (\ref{pivot1})
then (\ref{pivotbound}) is of order $n^{-4/5}$. Condition (\ref{fhatcond3}) now requires the bandwidth $h$ of the two density estimators to satisfy $n^{-4/5}=o(nh^7)$. This is achieved by choosing $h\gg n^{-9/35}$, thus allowing the rate optimal bandwidth   $h=c_2n^{-1/7}$.

\begin{rem}
The lower bounds for the integrated squared
error in deconvolution problems, as derived by Fan (1993), Theorem 2, also hold for uniform deconvolution.
Feuerverger et al. (2009) use this observation to show that their density estimator is rate optimal over Sobolev classes of densities.
If we compare our mean integrated squared error expansion with an optimal bandwidth of the form $h=cn^{-1/7}$
to the lower bound for a Sobolev class corresponding to twice differentiable densities then we
see that our estimator, for fixed $f$, also achieves the optimal rate $n^{-4/7}$.
\end{rem}

\section{A simulated example}\label{sim}

We use the average of $F^{-}_{nh}(x)$ and $F^{+}_{nh}(x)$ as
initial estimator in (\ref{finalest}) and (\ref{finalFest}).
Define
\begin{equation}\label{Fhat} \hat F_n(x) =F_{nh}^{(1/2)}=
\frac{1}{2}\,(F^{-}_{nh}(x) + F^{+}_{nh}(x)).
\end{equation}
The asymptotic variance of $\hat F_n(x)$ is equal to
$\frac{1}{4}\int w(u)^2du/nh$. Since the mean squared error
equals the sum of the squared bias and the variance we have
$$
\ex (\hat F_n(x)-F(x))^2= O(h^4)+O\Big(\frac{1}{nh}\Big),
$$
which asymptotically vanishes as long as $h\to 0$ and $nh\to \infty$.
If we choose a bandwidth of the form $h=c_1n^{-1/5}$ then the order of this mean
squared error is minimized. The minimal order is $n^{-4/5}$.
For the density estimators we choose a second bandwidth. We then have to ensure that
(\ref{fhatcond2}) holds, i.e. we should ensure $n^{-4/5}\ll nh^7$. This means that the bandwidth $h$ of
$f_{nh}(x)$ should satisfy $h\gg n^{-9/35}$. This is not an essential
restriction since the asymptotically optimal bandwidth that follows
from Theorem \ref{thm:1} is of order $n^{-1/7}$ and is thus
allowed.

As an illustration we have simulated a sample of size $n=500$ from
the convolution of the standard normal density ($f$) and the uniform
density.  The kernel function used is the biweight kernel
$$
w(x)=\frac{15}{16}\,(1-x^2)^2I_{[-1,1]}(x).
$$
The resulting estimates are given in Figures \ref{fig:1} and
\ref{fig:2}. For $f_{nh}^{-}$ and $f_{nh}^{+}$ we have used
the   bandwidth  $h=1$ and for $F_n$ we have chosen $h=0.7$. Indeed,
we see that the original estimates are relatively accurate in one
tail and almost periodic in the other tail. The combined estimate is
accurate  in both tails.

\begin{figure}[h]
$$
\includegraphics[width=60mm]{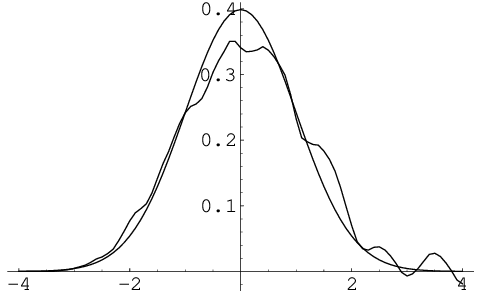} \quad
\includegraphics[width=60mm]{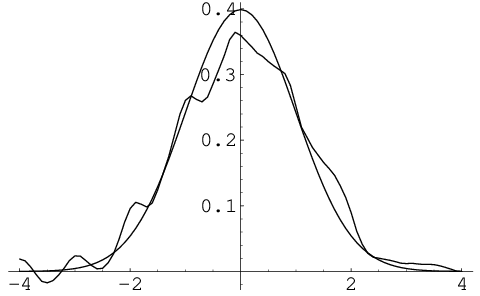}
$$
\caption[]{The estimates $f^{-}_{nh}$ and $f^{+}_{nh}$ and
the true density $f$, $h=1$. \label{fig:1}}
\end{figure}

\begin{figure}[h]
$$
\includegraphics[width=60mm]{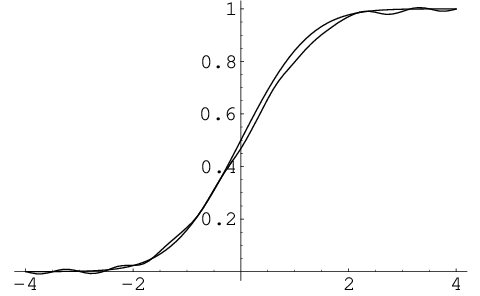}
\includegraphics[width=60mm]{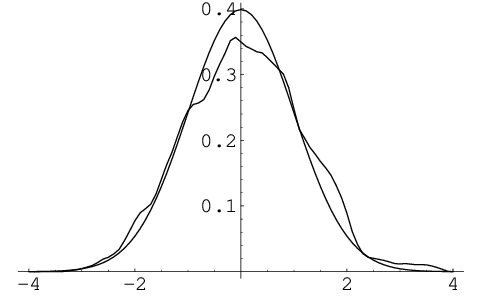}
$$
\caption[]{The estimate $F_n$, with $h=0.7$, and the final estimate
$f_{nh}$. \label{fig:2}}
\end{figure}

\bigskip

Next let us consider the estimator  of the distribution function. If we use   the specific estimator $\hat F_n(x)$ given by
 (\ref{Fhat}) then condition (\ref{fhatcond4})
requires $n^{-4/5}\ll nh^5$, which means $h\gg n^{-9/25}$. Again this is not an essential
restriction since the asymptotically optimal bandwidth that follows
from Theorem \ref{thm:2} is of order $n^{-1/5}$.

\begin{figure}[h]
$$
\includegraphics[width=60mm]{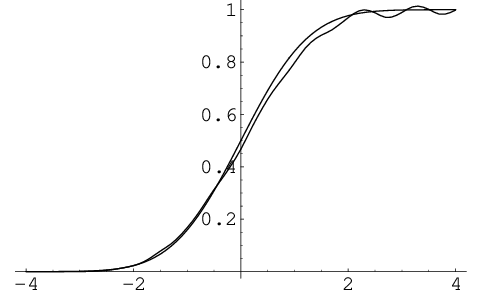}
\quad
\includegraphics[width=60mm]{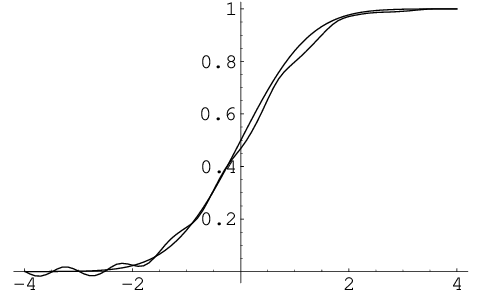}
$$
\caption[]{The estimates $F^{-}_{nh}$ and $F^{+}_{nh}$ and
the true distribution function $F$.
\label{fig:3}}
\end{figure}

\begin{figure}[h]
$$
\includegraphics[width=60mm]{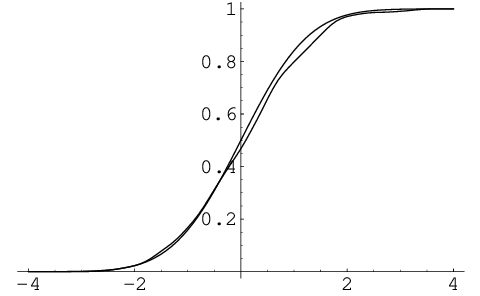}$$
\caption[]{The final estimate   $F_{nh}$ of $F$.
\label{fig:4}}
\end{figure}

Figures \ref{fig:3} and \ref{fig:4} give the estimates $F^{-}_{nh}$,
$F^{+}_{nh}$ and $F_{nh}$, based on the same sample of $n=500$ observations as above. Here
the bandwidth use is $h=0.7$. Again,
the original estimates are relatively accurate in one tail
and almost periodic in the other tail. The combined estimate is
accurate  in both tails.
Note the reduced variance in the tails of $F_{nh}$ compared to
that of $F_n$ in Figure \ref{fig:2}.

\section{Proofs}

\subsection{Proof of Theorem \ref{thm:1.1a} }

Note that
$$
f^{-}_{nh}(x)=\sum_{j=0}^\infty g'_{nh}(x-j)
=\frac{1}{nh^2}\sum_{i=1}^n\sum_{j=0}^\infty
w'\Big(\frac{x-j-X_i}{h}\Big)
$$
and
$$
f^{{+}}_{nh}(x)=-\sum_{j=1}^\infty g'_{nh}(x+j)
=-\frac{1}{nh^2}\sum_{i=1}^n\sum_{j=1}^\infty
w'\Big(\frac{x+j-X_i}{h}\Big).
$$
Write
$$
f^{(t)}_{nh}(x)= \sum_{i=1}^n \frac{1}{nh^2}\Big(t \sum_{j=0}^\infty
w'\Big(\frac{x-j-X_i}{h}\Big) -(1-t)\sum_{j=1}^\infty
w'\Big(\frac{x+j-X_i}{h}\Big)\Big) = \frac{1}{n}\sum_{i=1}^nU_{ih}(x)
$$
where
\begin{equation}\label{Uidef}
U_{ih}(x)=\frac{1}{h^2}\Big(t\sum_{j=0}^\infty
w'\Big(\frac{x-j-X_i}{h}\Big) -(1-t)\sum_{j=1}^\infty
w'\Big(\frac{x+j-X_i}{h}\Big)\Big).
\end{equation}

First we compute the expectations of the estimators. By
(\ref{inversion1}) we have
\begin{align}
\ex f^{-}_{nh}&(x)= \frac{1}{h^2}\sum_{j=0}^\infty \ex
w'\Big(\frac{x-j-X_1}{h}\Big) =\frac{1}{h^2}\sum_{j=0}^\infty
\int w'\Big(\frac{x-j-u}{h}\Big)g(u)du\nonumber\\
&=\frac{1}{h}\sum_{j=0}^\infty \int
w\Big(\frac{x-j-u}{h}\Big)g'(u)du =\sum_{j=0}^\infty
\frac{1}{h}\int w\Big(\frac{x-u}{h}\Big)g'(u-j)du\label{expleft}\\
&=\frac{1}{h}\int w\Big(\frac{x-v}{h}\Big)f(v)dv\nonumber
\end{align}
and similarly
\begin{equation}\label{expright}
\ex f^{+}_{nh}(x)=\frac{1}{h}\int
w\Big(\frac{x-u}{h}\Big)f(u)du.
\end{equation}
So the expectation of both $f^{-}_{nh}(x)$ and
$f^{+}_{nh}(x)$ is equal to the expectation of an ordinary
kernel estimator based on direct observations from $f$. From
(\ref{expleft}) and (\ref{expright}) we see that
\begin{equation}\label{expt}
\ex f^{(t)}_{nh}(x) =t \ex f^{-}_{nh}(x)+(1-t)\ex
f^{+}_{nh}(x) =\frac{1}{h}\int w\Big(\frac{x-u}{h}\Big)f(u)du.
\end{equation}
The bias expansion in the theorem now follows by standard arguments.

Similar to (\ref{expleft}) one can show $\ex U_{1h}(x)=O(1)$ if $f$ is
bounded on a neighborhood of $x$. The next lemma gives the even
moments of $U_{ih}(x)$.

\begin{lem}\label{momentslemma}
 For  $m $ even we have for $h\to 0$
\begin{equation}\label{moments}
\ex U_{ih}(x)^m=\frac{1}{h^{2m-1}}(t^m F(x) +(-1)^m(1-t)^m(1-F(x))\int
w'(v)^mdv+ O\Big( \frac{1}{h^{2m-2}}\Big).
\end{equation}
\end{lem}

\noindent{\bf Proof}

Note that
$$
  w'\Big(\frac{x-j_1-X_i}{h}\Big) w'\Big(\frac{x-j_2-X_i}{h}\Big)
=0
$$
if $j_1\not=j_2, j_1\in \mathbb Z, j_2\in \mathbb Z$ and $h<1/2$.
Similarly it is readily seen that the
 products of terms $w'\Big(\frac{x-j_l-X_i}{h}\Big)$
vanish if $h<1/2$ and if the $j_l$ are not all equal.

Now write

\begin{align*}
\ex&U_{ih}(x)^m=\frac{1}{h^{2m}}\ex \Big(t \sum_{j=0}^\infty
w'\Big(\frac{x-j-X_i}{h}\Big) -(1-t)\sum_{j=1}^\infty
w'\Big(\frac{x+j-X_i}{h}\Big)\Big)^m
\nonumber\\
&= \frac{1}{h^{2m}}\Big(t^m\sum_{j=0}^\infty  \ex
w'\Big(\frac{x-j-X_i}{h}\Big)^m + (-1)^m(1-t)^m\sum_{j=1}^\infty
\ex w'\Big(\frac{x+j-X_i}{h}\Big)^m\Big)\nonumber\\
&= \frac{1}{h^{2m-1}}\Big(t^m\sum_{j=0}^\infty \int
w'(v)^mg(x-j-hv)dv
\nonumber\\
& \quad\quad+ (-1)^m(1-t)^m\sum_{j=1}^\infty
\int w'(v)^mg(x+j-hv)dv\Big)\\
&=
\frac{1}{h^{2m-1}}\Big(t^m \int w'(v)^mF(x-hv)dv\\
&\quad\quad +
(-1)^m(1-t)^m\int w'(v)^m(1-F(x-hv))dv\Big)\nonumber\\
&= \frac{1}{h^{2m-1}}\,(t^m F(x) +(-1)^m(1-t)^m(1-F(x))\int
w'(v)^mdv+ O\Big(\frac{1}{h^{2m-2}}\Big). \nonumber
\end{align*}

\hfill$\Box$

\bigskip

For the variance of $f^{(t)}_{nh}(x)$ we get by Lemma
\ref{momentslemma}
\begin{align*}
\var f^{(t)}_{nh}&(x)= \frac{1}{n}\var(U_{1h}(x))
=\frac{1}{n}\Big(\ex U_{1h}(x)^2-(\ex U_{1h}(x))^2\Big)\nonumber\\
&\sim \frac{1}{nh^{3}}(t^2 F(x) +(1-t)^2(1-F(x))\int w'(v)^2dv.
\end{align*}

We will now check the Lyapunov condition for
$\frac{1}{n}\sum_{i=1}^n (U_{ih}(x)-\ex U_{ih}(x))$ to be asymptotically normal,
i.e. for some $\delta >0$ we have to check
$$
\frac{\ex |U_{1h}(x)-\ex
U_{1h}(x)|^{2+\delta}}{n^{\delta/2}(\var(U_{1h}(x)))^{1+\delta/2}}\to 0.
$$
Using $(a+b)^4\leq 2^3(a^4+b^4)$ we get, for suitable constants
$c_1$ and $c_2$ to be obtained from Lemma \ref{momentslemma},
$$
\frac{\ex (U_{1h}(x)-\ex U_{1h}(x))^4}{n(\var(U_{1h}(x)))^2} \leq \frac{2^3(\ex
U_{1h}(x)^4+(\ex U_{1h}(x))^4)}{n(\var(U_{1h}(x)))^2}\sim \frac{8c_1}{nhc_2^2}\to 0.
$$
This proves asymptotic normality of $(f^{(t)}_{nh}(x)-\ex
f^{(t)}_{nh}(x))/\sqrt{\var f^{(t)}_{nh}(x) }$ for fixed
$t$.\hfill$\Box$

\bigskip

\subsection{Proof of Theorem \ref{thm:1}}

We must show that we can replace $t=1-F(x)$ by a consistent
estimator. Write
\begin{equation}\label{decomp}
f_{nh}(x)=(1-\hat F_n(x)) f^{-}_{nh}(x) +  \hat
F_n(x)f^{+}_{nh}(x) =f_{nh}^{(1-F(x))}(x) +R_{nh}(x),
\end{equation}
where
\begin{equation}\label{RandWdefs}
R_{nh}(x)=(\hat F_n(x)-F(x)) S_{nh}(x)\quad\mbox{and}\quad
S_{nh}(x)=f^{+}_{nh}(x) -f^{-}_{nh}(x).
\end{equation}
Now write
$$
S_{nh}(x)=\frac{1}{n}\sum_{i=1}^n W_{ih}(x),
$$
where
$$
W_{ih}(x)=\frac{1}{h^2}\sum_{j=-\infty}^\infty
w'\Big(\frac{x-j-X_i}{h}\Big).
$$
The next lemma establishes some properties of
$S_{nh}(x)$.

\begin{lem}\label{Wlemma}
We have $\ex S_{nh}(x)=0$, $\ex W_{ih}(x)^m=\frac{1}{h^{2m-1}}\int w'(u)^mdu$ and
$$
\sqrt{nh^3}S_{nh}(x)\convd
N\Big(0,\int w'(u)^2du\Big).
$$
The distributions of the random variables $W_{ih}(x)$ and $S_{nh}(x)$ are independent of $x$.
\end{lem}
\noindent{\bf Proof}

The first statement follows from (\ref{expleft}) and (\ref{expright}).
The second statement follows from a computation similar to the one in the proof of
Lemma \ref{momentslemma}. Asymptotic normality can be proved as in the proof of
Theorem \ref{thm:1.1a}.

The fact that the distribution does not depend on $x$ can be seen by writing
$$
\sum_{j=-\infty}^\infty
w'\Big(\frac{x-j-X_i}{h}\Big)=\sum_{j=-\infty}^\infty
w'\Big(\frac{x-j-Y_i-Z_i}{h}\Big).
$$
Given $x$ and $Y_i$ this sum equals a periodic function
with period one evaluated at $Z_i$. Since $Z_i$ is Uniform$[0,1)$ distributed its distribution
does not depend on $x$ and $Y_1$.

\hfill$\Box$

\bigskip

By (\ref{fhatcond}) we
have $\hat F_n(x)-F(x)\convp 0$ and hence by Slutsky's theorem $
\sqrt{nh^3}R_{nh}(x)\convp 0. $ Furthermore by the Cauchy-Schwarz
inequality
$$
\ex \sqrt{nh^3}|R_{nh}(x)| \leq \sqrt{nh^3} (\ex (\hat
F_n(x)-F(x))^2)^{1/2} (\ex (S_{nh}(x))^2)^{1/2}\to 0.
$$
This shows that $\sqrt{nh^3}(f_{nh}(x)-\ex f_{nh}(x))$ has the
same asymptotic normal distribution as
$\sqrt{nh^3}(f_{nh}^{(1-F(x))}(x)-\ex f_{nh}^{(1-F(x))}(x))$,
which proves the first statement of the theorem.

To prove the second statement note that by (\ref{expt}) and
 a standard argument in kernel estimation we have
\begin{equation}\label{expft}
\ex f_{nh}^{(1-F(x))}(x) = \frac{1}{h}\int
w\Big(\frac{x-u}{h}\Big)f(u)du = f(x) + \frac{1}{2}h^2f''(x)\int
v^2w(v)dv+ o(h^2).
\end{equation}
Furthermore
\begin{align}
\ex |R_{nh}(x)| &\leq  (\ex (\hat F_n(x)-F(x))^2)^{1/2}
(\ex (S_{nh}(x))^2)^{1/2} \nonumber\\
&= o\Big(\sqrt{nh^7}\Big)O\Big(\frac{1}{\sqrt{nh^3}}\Big)=o(h^2).
\label{biasest}\end{align} Together (\ref{expft}) and
(\ref{biasest}) prove the second statement of the theorem.
\hfill$\Box$

\subsection{Proof of Theorem \ref{thm:2.1b}}

We copy the proof of Theorem \ref{thm:1}. Note that
$$
F^{-}_{nh}(x)=\sum_{j=0}^\infty g_{nh}(x-j)
=\frac{1}{nh}\sum_{i=1}^n\sum_{j=0}^\infty
w\Big(\frac{x-j-X_i}{h}\Big)
$$
and
$$
F^{+}_{nh}(x)=1-\sum_{j=1}^\infty g_{nh}(x+j)
=1-\frac{1}{nh}\sum_{i=1}^n\sum_{j=1}^\infty w
\Big(\frac{x+j-X_i}{h}\Big).
$$
Write
\begin{align*}
F^{(t)}_{nh}&(x)=t F^{-}_{nh}(x) +(1-t) F^{+}_{nh}(x)
\nonumber\\
&= \sum_{i=1}^n \frac{1}{nh}\Big(t \sum_{j=0}^\infty
w\Big(\frac{x-j-X_i}{h}\Big)
-(1-t)\sum_{j=1}^\infty w\Big(\frac{x+j-X_i}{h}\Big)\Big)+1-t\\
&= \frac{1}{n}\sum_{i=1}^nV_{ih}(x)+1-t \nonumber
\end{align*}
where
$$
V_{ih}(x)=\frac{1}{h}\Big(t \sum_{j=0}^\infty w\Big(\frac{x-j-X_i}{h}\Big)
-(1-t)\sum_{j=1}^\infty w\Big(\frac{x+j-X_i}{h}\Big)\Big).
$$

First we compute the expectations of the estimators. By
(\ref{inversion1}) we have
\begin{align}
\ex F^{-}_{nh}&(x)= \frac{1}{h}\sum_{j=0}^\infty \ex
w\Big(\frac{x-j-X_1}{h}\Big) =\frac{1}{h}\sum_{j=0}^\infty
\int w\Big(\frac{x-j-u}{h}\Big)g(u)du\nonumber\\
&= \sum_{j=0}^\infty
\frac{1}{h}\int w\Big(\frac{x-u}{h}\Big)g(u-j)du\label{Fexpleft}\\
&=\frac{1}{h}\int w\Big(\frac{x-u}{h}\Big)F(u)du\nonumber
\end{align}
and similarly
\begin{equation}\label{Fexpright}
\ex F^{+}_{nh}(x)=\frac{1}{h}\int
w\Big(\frac{x-u}{h}\Big)F(u)du.
\end{equation}
 Since it is a convex combination of $F^{-}_{nh}(x)$ and $F^{+}_{nh}(x)$
the expectation of $F^{(t)}_{nh}(x)$ is also equal to
(\ref{Fexpleft}) and (\ref{Fexpright}).

The equivalent to Lemma \ref{momentslemma} for $V_{ih}(x)$ is
\begin{equation}\label{moments2}
\ex V_{ih}(x)^m=\frac{1}{h^{m-1}}(t^m F(x) +(-1)^m(1-t)^m(1-F(x))\int
w(v)^mdv+ o\Big( \frac{1}{h^{m-1}}\Big)
\end{equation}
which follows by replacing $w'$ by $w$ and replacing $1/h^2$ by
$1/h$ in the proof.

For the variance of $F^{(t)}_{nh}(x)$ we then get
\begin{align*}
\var F^{(t)}_{nh}&(x)= \frac{1}{n}\var(V_{1h}(x))
=\frac{1}{n}\Big(\ex V_{1h}(x)^2-(\ex V_{1h}(x))^2\Big)\nonumber\\
&\sim \frac{1}{nh}(t^2 F(x) +(1-t)^2(1-F(x))\int w(v)^2dv.
\end{align*}

The Lyapunov condition for $\frac{1}{n}\sum_{i=1}^n (V_{ih}(x)-\ex V_{ih}(x))$
to be asymptotically normal can be checked as in the proof of
Theorem \ref{thm:1}. This proves asymptotic normality of
$(F^{(t)}_{nh}(x)-\ex F^{(t)}_{nh}(x))/\sqrt{\var F^{(t)}_{nh}(x) }$
for fixed $t$. \hfill$\Box$

\subsection{Proof of Theorem \ref{misethm}}

Recall that by (\ref{decomp}) we have
$$
f_{nh}(x)=f_{nh}^{(1-F(x))}(x) +R_{nh}(x),
$$
where
$$
R_{nh}(x)=(\hat F_n(x)-F(x)) S_{nh}(x)\quad\mbox{and}\quad
S_{nh}(x)=f^{+}_{nh}(x) -f^{-}_{nh}(x).
$$
We decompose the mean integrated squared error as follows
\begin{align}
\MISE&(h) =\int \ex \Big(f_{nh}^{(1-F(x))}(x)-f(x) + R_{nh}(x)\Big)^2dx \nonumber\\
&= \int \ex \Big(f_{nh}^{(1-F(x))}(x)-f(x)\Big)^2dx +\int \ex R_{nh}(x)^2dx\label{ter1}\\
&\quad\quad+2\int \ex \Big((f_{nh}^{(1-F(x))}(x)-f(x)\Big) R_{nh}(x)) dx\label{ter2}.
\end{align}
The mean integrated squared error of $f_{nh}^{(1-F(x))}$ can be written as integrated
squared bias plus integrated squared variance. We have
 $$
\int \ex (f_{nh}^{(1-F(x))}(x)-f(x))^2dx =
\int (\ex f_{nh}^{(1-F(x))}(x)- f(x))^2 dx+ \int \var f_{nh}^{(1-F(x))}(x)dx.
$$

We have already noted in (\ref{expt}) that the expectation of $f_{nh}^{(1-F(x))}(x)$
is equal to the expectation of a standard
kernel estimator.
By Theorem 2.1.7 of Prakasa Rao (1983), or the original proof of Nadaraya,
 we have the standard expansion for integrated squared
bias of $f_{nh}^{(1-F(x))}(x)$, i.e.
$$
\int (\ex f_{nh}^{(1-F(x))}(x)-f(x))^2dx
=\frac{1}{4}h^4\int f''(x)^2dx\Big(\int v^2w(v)dv\Big)^2+
o(h^4).
$$
Next consider the integrated variance. We have, with $U_{ih}(x)$ as in
(\ref{Uidef}),
\begin{equation}\label{var2}
\var f_{nh}^{(1-F(x))}(x)=
 \frac{1}{n}\var(U_{1h}(x))
=\frac{1}{n}\Big(\ex U_{1h}(x)^2-(\ex U_{1h}(x))^2\Big),
\end{equation}
where by the proof of Lemma \ref{momentslemma}, with $t=1-F(x)$,
$$
\ex U_{1h}(x)^m=\frac{1}{h^{2m-1}}\Big((1-F(x))^m \int w'(v)^mF(x-hv)dv
   +
(-1)^mF(x)^m\int w'(v)^m(1-F(x-hv))dv\Big)
$$
Now use
$$
F(x-hv)= F(x) - hv   \int_0^1 f(x-thv)dt
$$
to get
\begin{align*}
\ex U_{1h}(x)^2=&\frac{1}{h^{3}}F(x)(1-F(x)) \int w'(v)^2dv\\
 -\frac{1}{h^{2}}&((1-F(x))^2+F(x)^2 )\int_{-1}^1\int_0^1 vw'(v)^2f(x-thv)dvdt.
\end{align*}
The integral    with respect to $x$ of the first
term is finite by the condition $\ex Y<\infty$.
The integral    with respect to $x$ of the second
term is finite by the fact that $|(1-F(x))^2+F(x)^2|$ is bounded by two and
Fubini's theorem. Similarly, for the term $\ex U_{1h}(x)$ in (\ref{var2})
we get by the Cauchy Schwartz inequality and Fubini's theorem
\begin{align*}
\int (\ex U_{1h}(x))^2&dx = \int\Big( \int_{-1}^1\int_0^1  vw'(v)f(x-thv)dvdt\Big)^2dx\\
&\leq \int\Big( \int_{-1}^1\int_0^1  v^2w'(v)^2dvdt\int_{-1}^1\int_0^1f(x-thv)^2dvdt\Big) dx\\
&=\int_{-1}^1\int_0^1  v^2w'(v)^2dvdt\int\int_{-1}^1\int_0^1f(x-thv)^2dxdvdt\\
&=2\int_{-1}^1 v^2w'(v)^2dv \int f(x)^2dx.
\end{align*}

Finally this gives
$$
\int \var f_{nh}^{(1-F(x))}(x)dx=\frac{1}{nh^3}F(x)(1-F(x)) \int w'(v)^2dv
+ O\Big(\frac{1}{nh^2}\Big).
$$

For the integrated expected squared remainder term in (\ref{ter1}) we have
\begin{align*}
\int \ex R_{nh}&(x)^2dx=\int \ex (\hat F_n(x)-F(x))^2  S_{nh}(x)^2dx\\
&\leq\int (\ex (\hat F_n(x)-F(x))^4)^{1/2} (\ex S_{nh}(x)^4)^{1/2}dx\\
&=  o(nh^7)O\Big(\frac{1}{n h^3}\Big)=o(h^4),
\end{align*}
since by Lemma \ref{Wlemma}
\begin{align*}
(\ex S_{nh}(x)^4&)^{1/2}=\Big\{\frac{1}{n^4}\Big[n\ex W_{1h}(x)^4+ 3n(n-1)(\ex W_{1h}(x))^2\Big]\Big\}^{1/2}\\
&=\Big\{\frac{1}{n^4}\Big[n\Big(\frac{1}{h^7}\int w'(v)^4dv\Big)+ 3n(n-1)\Big(\frac{1}{h^3}\int w'(v)^2dv\Big)^2\Big]
\Big\}^{1/2}=O\Big(\frac{1}{nh^3}\Big).
\end{align*}
The proof of the theorem is completed by noting that the cross product term
  (\ref{ter2}) is negligible with respect to the first term  (\ref{ter1}) by
the Cauchy Schwartz inequality.
\hfill$\Box$

\subsection{Proof of Lemma \ref{pivotok}}

Write
$$
F_{nh}^{(H(x))}(x)-F(x)=F_{nh}^{(H(x))}(x)-\ex F_{nh}^{(H(x))}(x) +\ex F_{nh}^{(H(x))}(x)-F(x)
$$
By the triangle inequality we have
$$
\Big(\ex(F_{nh}^{(H(x))}(x)-F(x))^4\Big)^{1/4}\leq   \Big(\ex(F_{nh}^{(H(x))}(x)-\ex F_{nh}^{(H(x))}(x))^4\Big)^{1/4}+\Big(\ex F_{nh}^{(H(x))}(x)-F(x)\Big) .
$$
So by $(a+b)^2\leq 2(a^2+b^2), a,b\geq 0$,  we also have
$$
\Big(\ex(F_{nh}^{(H(x))}(x)-F(x))^4\Big)^{1/2}\leq 2  \Big(\ex(F_{nh}^{(H(x))}(x)-\ex F_{nh}^{(H(x))}(x))^4\Big)^{1/2}
+2 \Big(\ex F_{nh}^{(H(x))}(x)-F(x)) \Big)^{2}.
$$

Hence it suffices to prove the bound of the lemma for the fourth power of the error and the usual square of the bias separately.

In the proof of Theorem \ref{thm:2.1b} we have seen that  $F_{nh}^{(H(x))}(x)$ can be written as
$$
 F^{(H(x))}_{nh}(x)=\frac{1}{n}\sum_{i=1}^nV_{ih}(x)+H(x)
$$
where
$$
V_{ih}(x)=\frac{1}{h}\Big((1-H(x))\sum_{j=0}^\infty w\Big(\frac{x-j-X_i}{h}\Big)
-H(x)\sum_{j=1}^\infty w\Big(\frac{x+j-X_i}{h}\Big)\Big).
$$
and that we have
$$
\ex F^{(H(x))}_{nh}(x)=\frac{1}{h}\int
w\Big(\frac{x-u}{h}\Big)F(u)du.
$$
Following the same arguments as in the proof of Theorem 2.1.7, the MISE expansion for kernel estimators,
 of Prakasa Rao (1983) we have
\begin{equation}
\int_{-\infty}^\infty (\ex F_{nh}^{(H(x))}(x)-F(x))^2 dx= \frac{1}{4} h^4 \Big(\int_{-\infty}^\infty f'(x)^2dx\Big)\Big(\int_{-\infty}^\infty  v^2w(v)dv\Big)^2 +
o(h^4).
\end{equation}
In order to deal with the error part we wite
$$
F^{(H(x))}_{nh}(x)-\ex F^{(H(x))}_{nh}(x)=\frac{1}{n}\sum_{i=1}^n \tilde V_{ih}(x),
$$
where $\tilde V_{ih}(x) = V_{ih}(x)-\ex V_{ih}(x)$.  Since   $\ex \tilde V_{ih}(x)$ equals zero we have
$$
\ex\Big(\frac{1}{n}\sum_{i=1}^n \tilde V_{ih}(x)\Big)^4=\frac{1}{n^3}\ex\Big(\tilde V_{1h}(x)^4\Big)  +
3\,\frac{n-1}{n^3} \Big(\ex\Big(\tilde V_{1h}(x)^2\Big) \Big)^2.
$$
From (\ref{moments2}) we get
$$
\frac{1}{n^3}\ex\Big(\tilde V_{1h}(x)^4\Big) \sim \frac{1}{n^3}\ex\Big(V_{1h}(x)^4\Big) \sim \frac{c_1}{n^3h^3}\,\Big( (1-H(x))^4F(x)+ H(x)^4(1-F(x))\Big)
$$
and
$$
3\,\frac{n-1}{n^3} \Big(\ex\Big(\tilde V_{1h}(x)^2\Big) \Big)^2=
3\,\frac{n-1}{n^3}\big(\var (V_{1h}(x))^2\Big) \sim \frac{c_2}{n^2h^2}\,\Big( (1-H(x))^2F(x)+ H(x)^2(1-F(x))\Big)^2 ,
$$
for certain constants $c_1$ and $c_2$.
Since the square roots of the functions on the right hand side are integrable we now have
$$
\int_{\infty}^\infty \Big(\ex(F^{(H(x))}_{nh}(x)-\ex F^{(H(x))}_{nh}(x))^4\Big)^{1/2}dx=O\Big(\frac{1}{nh}\Big)
$$

Summarizing we get
$$
\int_{\infty}^\infty  \Big(\ex(F_{nh}^{(H(x))}(x)-F(x))^4\Big)^{1/2}dx = O(h^4) + O\Big(\frac{1}{nh}\Big),
$$
which completes the proof.
\hfill$\Box$

\section{The limit variance of the smoothed NPMLE}\label{GJvar}

In this section we will show that the limit variance of the smoothed NPMLE is equal to the
limit variance of our optimally combined kernel estimator.

Let the distribution induced by $F$ have support $[0,M)$ for some $M>0$ and let
$m$ denote the largest integer strictly smaller than $M+1$.
The asymptotic variance of the smoothed NPMLE in Theorem 2 in Groeneboom and Jongbloed (2003) is, in their notation where 
$t$ stands for our $x$ in Theorem \ref{thm:1}, defined as
\begin{equation}\label{var1}
\sigma^2=\lim_{h\downarrow 0}
\int \theta^2_{h,t,F}dG,
\end{equation}
with the function $\theta_{h,t,F}$, for $0\leq t<M$, defined by 
\begin{equation}
\theta_{h,t,F}(x+k)=\left\{
\begin{array}{ll}
  \sum_{i=0}^m(1-F(x+i))w_h'(t-(x+i)) & ,\ \mbox{if}\ x\in[0,1], \ k=0, \\
 -\sum_{i=0}^{k-1}w_h'(t-(x+i))+\theta_{h,t,F}(x) &,\ \mbox{if}\ x\in[0,1], \ k=1,\ldots,m,
\end{array}
\right.
\end{equation}
where $w_h(\cdot)=w(\cdot/h)/h$.

\begin{lem}
The asymptotic variance (\ref{var1}) is equal to $F(t)(1-F(t))\int w'(u)^2du.$
\end{lem}

\noindent{\bf Proof}
We write the integral in (\ref{var1}) as
\begin{align*}
\int &\theta^2_{h,t,F}dG  = \sum_{k=0}^{m}\int_0^1\theta^2_{h,t,F}(x+k)g(x+k)dx\\
&=\int_0^1 (F(x)-F(x-1))\Big[(1-F(x))w_h'(t-x)+(1-F(x+1))w_h'(t-x-1)\\
&\quad\quad+(1-F(x+2))w_h'(t-x-2)+\ldots+(1-F(x+m))w_h'(t-x-m)\big]^2dx\\
&+\int_0^1 (F(x+1)-F(x))\Big[-F(x)w_h'(t-x)+(1-F(x+1))w_h'(t-x-1)\\
&\quad\quad+(1-F(x+2))w_h'(t-x-2)+\ldots+(1-F(x+m))w_h'(t-x-m)\big]^2dx\\
&+\int_0^1 (F(x+2)-F(x+1))\Big[-F(x)w_h'(t-x) -F(x+1) w_h'(t-x-1)\\
&\quad\quad+(1-F(x+2))w_h'(t-x-2)+\ldots+(1-F(x+m))w_h'(t-x-m)\big]^2dx\\
&+ \ldots \\
&+\int_0^1 (F(x+m)-F(x+m-1))\Big[-F(x)w_h'(t-x) -F(x+1) w_h'(t-x-1)\\
&\quad\quad -F(x+2) w_h'(t-x-2)+\ldots -F(x+m) w_h'(t-x-m)\Big]^2dx.
\end{align*}
The next step is to write out the squares, which we leave to the reader.
Let $l\leq t<l+1$ for some integer $l$, then, since $0\leq t-l <1$ and $x\in[0,1]$, only the terms containing $w_h'(t-x-l)^2$
will yield a non zero contribution for $h$ small enough.
This contribution is, for $h$ to zero, equal to
\begin{align*}
&\int_0^1 \sum_{j=0}^l(F(x+j)-F(x+j-1))(1-F(x+l))^2w_h'(t-x-l)^2dx \\
&\quad\quad+\int_0^1 \sum_{j=l+1}^ {m } (F(x+j)-F(x+j-1))F(x+l)^2w_h'(t-x-l)^2dx\\
&=\int_0^1 F(x+l)(1-F(x+l))^2w_h'(t-x-l)^2dx\\
&\quad\quad+\int_0^1  (1-F(x+l ))F(x+l)^2w_h'(t-x-l)^2dx\\
&\sim \frac{1}{h^3}\, \Big(F(t)(1-F(t))^2+(1-F(t))F(t)^2\Big) \int w'(u)^2du\\
&=\frac{1}{h^3}\,  F(t)(1-F(t)) \int w'(u)^2du.
\end{align*}
Here we have used an expansion of the integral which is standard in kernel estimation theory.

Taking the limit for $h$ to zero as in (\ref{var1}) now yields the result.
\hfill$\Box$

\end{document}